\documentclass[letterpaper, 10 pt, conference]{ieeeconf}  

\IEEEoverridecommandlockouts                              
\usepackage{amsmath, amssymb, amsxtra, mathrsfs}
\usepackage{graphics, epsfig, here, epstopdf, fancyhdr,graphicx}
\usepackage{bbm}
\usepackage{color}
\usepackage{multirow}
\usepackage[all]{xy}

\newcommand{\tS}{\textstyle}

\newcommand{\pf}{\noindent \textit{Proof:}\quad}

\newcommand{\bN}{\mathbb{N}}

\newcommand{\bR}{\mathbb{R}}

\newcommand{\cC}{\mathcal{C}}
\newcommand{\cE}{\mathcal{E}}

\newcommand{\cL}{\mathcal{L}}

\newcommand{\cV}{\mathcal{V}}

\newcommand{\cM}{\mathcal{M}}
\newcommand{\cD}{\mathcal{D}}

\newcommand{\cT}{\mathcal{T}}

\newcommand{\fA}{\,\forall\,}

\newcommand{\dS}{\displaystyle}

\newcommand{\bp}[1]{{\left(#1\right)}}

\newcommand{\rM}[1]{{\mathrm{#1}}}

\newcommand{\mB}[1]{{\mathbf{#1}}}

\newcommand{\Bone}{\mB{1}}

\newcommand{\vphi}{\varphi}

\newcommand{\cG}{\mathcal{G}}

\newcommand{\pp}{\partial}

\newcommand{\wh}[1]{{\widehat{#1}}}

\newcommand{\thalf}{\tfrac{1}{2}}

\newcommand{\rd}{\mathrm{d}}
\newcommand{\ddx}[1]{{\frac{\rd }{\rd #1}}}

\usepackage{algorithmic}
\usepackage{epstopdf}

\title{\LARGE \bf
Optimal Control of Transient Flow in Natural Gas Networks
}

\author{Anatoly Zlotnik$^\dagger$ and Michael Chertkov$^\ddag$ and Scott Backhaus$^\S$
\thanks{$^\dag$azlotnik@lanl.gov, \,\, Center for Nonlinear Studies, Theoretical Division, Los Alamos National Laboratory, Los Alamos, NM 87545}
\thanks{$^\ddag$chertkov@lanl.gov, \,\ Physics of Condensed Matter \& Complex Systems, Theoretical Division, Los Alamos National Laboratory, Los Alamos, NM 87545}
\thanks{$^\S$backhaus@lanl.gov, \,\, Condensed Matter \& Magnet Science,
    Materials Physics and Applications Division, Los Alamos National Laboratory, Los Alamos, NM 87545}}

\newtheorem{thm}{Theorem}

\newtheorem{prop}[thm]{Proposition}

\newtheorem{rem}{Remark}

\begin{document}

\maketitle
\thispagestyle{empty}
\pagestyle{empty}

\begin{abstract}
We outline a new control system model for the distributed dynamics of compressible gas flow through large-scale pipeline networks with time-varying injections, withdrawals, and control actions of compressors and regulators.  The gas dynamics PDE equations over the pipelines, together with boundary conditions at junctions, are reduced using lumped elements to a sparse nonlinear ODE system expressed in vector-matrix form using graph theoretic notation. This system, which we call the reduced network flow (RNF) model, is a consistent discretization of the PDE equations for gas flow.  The RNF forms the dynamic constraints for optimal control problems for pipeline systems with known time-varying withdrawals and injections and gas pressure limits throughout the network.  The objectives include economic transient compression (ETC) and minimum load shedding (MLS), which involve minimizing compression costs or, if that is infeasible, minimizing the unfulfilled deliveries, respectively.  These continuous functional optimization problems are approximated using the Legendre-Gauss-Lobatto (LGL) pseudospectral collocation scheme to yield a family of nonlinear programs, whose solutions approach the optima with finer discretization.  Simulation and optimization of time-varying scenarios on an example natural gas transmission network demonstrate the gains in security and efficiency over methods that assume steady-state behavior.
\end{abstract}

\section{Introduction} \label{secintro}

New emissions restrictions and the resulting push towards cleaner electric power sources, as well as increased supplies of natural gas in the United States, have compelled the installation of gas-fired electric power plants for the vast majority of new generating capacity over the past 15 years \cite{lyons13}.  Such generators can quickly adjust their output, hence they are often dispatched to balance out the fluctuating and highly variable production of uncontrollable renewable energy sources such as wind and solar \cite{li08,mitei14}.  This growing and increasingly time-variable natural gas consumption creates a significant impact on the pressure and flow throughout associated natural gas transmission networks.  These conditions contrast with historically slower and smaller variations in withdrawals by local distribution companies (LDCs), which allowed gas pipeline operators to assume, for day-ahead planning purposes, that consumption remains constant throughout the day.  Today's new complexities cause interactions on previously distinct spatiotemporal scales, which present risks and disruptive challenges that invalidate many traditional approaches for design, risk assessment, and operation of natural gas transmission systems \cite{chertkov14,tabors14,chertkov15hicss}.

Historically, natural gas was predominantly withdrawn from transmission systems by LDCs and industrial consumers in a predictable way and with relatively little variation over a day.  It was traded using day-ahead contracts for fixed deliveries with the assumption that injections and withdrawals will remain nearly constant.  Accordingly, early studies \cite{wong68,rothfarb70,luongo91} focused on optimizing steady-state gas flows, for which the state equations are algebraic relations.  Recent efforts have improved and scaled up optimization techniques for similar problems  \cite{midthun09,borraz10phd,misra15,babonneau14}. However, the steady-state assumption no longer represents realistic operating conditions because of the economic and technological trends discussed above \cite{peters12}.  Intermittent dispatch of gas-fired power plants
creates rapidly changing gas withdrawals from transmission pipelines, and results in
stresses that are increasingly difficult to contain within system design limits using current ad-hoc methods.  The growing reliance on natural gas for electricity production has led to strong coupling between electric power and natural gas infrastructures, and has created a need for secure gas transmission to prevent electric load shedding because of interrupted gas deliveries to generators \cite{mitei14}.  This new challenge compels our investigation of new techniques for modeling and optimal control of dynamic flows of compressible gas in large-scale pipeline networks.

The transient flow of natural gas in a transmission pipeline can be represented by simplifications to the Euler equations for compressible gas flow in one-dimension \cite{wylie78,osiadacz84}.  For physically relevant pipeline parameters, this PDE system is defined over very large scales in both distance and time, and is highly nonlinear even with numerous physical modeling simplifications \cite{dorao11}.  Transient flows in pipelines on the scale of thousands of miles are thus problematic to simulate, and many methods have been proposed \cite{thorley87,grundel13a,herty10}. 
Indeed, pipeline simulation is an area of active research of interest to gas system operators \cite{chapman05}.

The difficulty of characterizing gas network dynamics presents challenges for engineering, design, and operation of natural gas transmission systems under transient conditions.  In the related optimization problems, the PDEs representing the dynamics are incorporated as constraints that must be satisfied over widely distributed space and time domains \cite{moritz07phd}, and their nonlinearity makes computational tractability a challenge. Previous studies examined optimization of multi-day operations of gas pipeline networks involving transient flows \cite{rachford00,ehrhardt05,steinbach07pde}, including work on economic model predictive control \cite{gopalakrishnan13,devine14}.  In these studies the PDE constraints for gas flow are represented using implicit first-order schemes in space and time, which result in very large-scale problems because of the fine discretization required to adequately resolve transients on the time-scales of interest.

In this manuscript, we develop a modeling and control framework that closely represents the physical phenomena in gas pipeline networks.  We also present methods for simulation and optimization of dynamic compressible gas flows over such systems using nodal actuators, which provides unprecedented gains in efficiency and scalability over the previous work listed above.  For example, the optimal gas flow (OGF) \cite{misra15} determines optimal compressor station set-points that balance constant injections and withdrawals over a network while satisfying system pressure constraints. Here, we extend this concept to the transient case.  We examine the objectives of economic transient compression (ETC) and minimum load shedding (MLS), which involve minimizing compression costs or unfulfilled deliveries to customers with non-firm contracts, respectively.  In the transient regime, these problems are formulated as PDE-constrained optimal control problems (OCPs).   The PDE constraints are approximated by a new control system model, the reduced network flow (RNF), derived from a model reduction of gas network dynamics \cite{grundel13a}.  The approximation used for the RNF has been validated through comparison with an operator split-step numerical solution to the one-dimensional PDE system for a single pipe with transient boundary conditions \cite{zlotnik15dscc}.  The modified OCP is then approximated with a nonlinear program (NLP) using pseudospectral discretization \cite{ruths11cdc}.  The decision variables are coefficients of a polynomial expansion that approximates the solution to the OCP.  Our approach provides several advantages relative to previous methods.  The representation of continuous dynamics using polynomials gives spectral accuracy, yielding comparable fidelity using coarser discretization with far fewer decision variables.  The RNF equations can also be integrated with an ODE solver to validate solutions to the NLP.

The manuscript is organized as follows.  In Section \ref{sec:pipeline}, we summarize the physics of compressible gas flow in pipelines.  Section \ref{sec:netflow} contains an optimal control formulation for dynamic flows on networks subject to time-varying injections, withdrawals, and actions of compressors and regulators, and defines the ETC and MLS OCPs.  In Section \ref{sec:model}, we derive the RNF control system model, and provide consistency results.  In Section \ref{sec:pseudospectral}, we summarize the Legendre-Gauss-Lobatto (LGL) pseudospectral collocation method for optimal control.  Section \ref{sec:implement} describes implementation of the LGL scheme to approximate the reduced OCPs in Section \ref{sec:model} with NLPs.  Section \ref{sec:examples} contains a case study of the ETC and MLS problems for an example gas network.  Section \ref{sec:conc} contains a brief discussion of the results and their implications.


\vspace{-1ex}
\section{Gas Pipeline Dynamics} \label{sec:pipeline}

The dissipative flow of a compressible gas in a horizontal pipeline with slow transients that do not excite waves or shocks is adequately described by a simplification of the Euler equations in one dimension \cite{thorley87,herty10,osiadacz84}, given by
\begin{align}
 \dS \pp_t\rho+\pp_x\vphi &= 0,  \label{euler1a} \\
\dS\pp_t\vphi + a^2\pp_x \rho &= -\frac{\lambda}{2D}\frac{\vphi |\vphi|}{\rho}. \label{euler1b}
\end{align}
The variables are mass flux $\vphi$ and density $\rho$, defined on a domain $x\in [0,L]$ at time $t$.   The parameters are the friction factor $\lambda$, pipe diameter $D$, and speed of sound $a$.  The term on the right hand side of \eqref{euler1b} aggregates friction effects.  We have assumed that gas pressure $p$ and density $\rho$ satisfy the relation $p=a^2\rho$ with $a^2=ZRT$, where $Z$, $R$, and $T$,  are the gas compressibility factor, ideal gas constant, and constant temperature, respectively.   Equation \eqref{euler1b} is valid in the regime when changes in the boundary conditions are sufficiently slow to not excite propagation of sound waves.
The gas dynamics on a pipeline segment are represented using \eqref{euler1a}-\eqref{euler1b}, with a unique solution when any two of the boundary conditions $\rho(t,0)=\rho_{0}(t)$, $\vphi(t,0)=\vphi_{0}(t)$, $\rho(t,L)=\rho_{L}(t)$, or $\vphi(t,L)=\vphi_{L}(t)$ are specified.


Because of friction, the pressure of gas flowing through a pipeline gradually decreases, and is boosted by compressors so that it exceeds the minimum for delivery to customers.  Compressor stations are controllable actuators used to manipulate the state of the transmission system.    The physical size of a station is very small compared to a pipeline, hence we model compressor/regulator action as conservation of flow and a multiplicative change in density at a point $x=c$.  Specifically, $\rho(t,c^+)=\alpha(t)\rho(t,c^-)$ and $\vphi(t,c^+)=\vphi(t,c^-)$, where $\alpha(t)$ is a time-dependent compression factor.  We denote $h(c^-)=\lim_{x\nearrow c}h(x)$ and $h(c^+)=\lim_{x\searrow c}h(x)$.
The compression power is proportional to
\begin{align} \label{obj1}
C\propto \eta^{-1}|\vphi(t,c)|(\max\{\alpha(t),1\}^{2m}-1)
\end{align}
with $0<m<(\gamma-1)/\gamma<1$ where $\gamma$ is the heat capacity ratio and $\eta$ is the compressor efficiency \cite{wong68,misra15}.

\section{Network Flow Control Formulation} \label{sec:netflow}

We consider a network of pipelines that are connected at junctions where gas flow can be compressed, or gas can be withdrawn from or injected into the system.  Between junctions, the mass flux and density evolve according to \eqref{euler1a}-\eqref{euler1b}. This collection of segments connected at junctions is considered as a directed graph $\cG=(\cV,\cE)$, where each segment is an edge $\{i,j\}\in\cE$ in the set of edges $\cE$ that connects junctions $i,j\in\cV$ in the set of nodes $\cV$.  The instantaneous state within the edge $\{i,j\}$ is characterized by the density $\rho_{ij}$ and flux $\vphi_{ij}$ defined on a time interval $\cT=[0, T]$ and on the distance variable $x_{ij} \in[0, L_{ij}]=\cL_{ij}$, where $L_{ij}$ is the length of edge $\{i,j\}$.  Defining the domain of the PDE solution for the $\{i,j\}\in\cE$ as $\cD_{ij}=\cT\times \cL_{ij}$, the state functions are expressed as $\rho_{ij}:\cD_{ij}\to\bR_+$ and $\vphi_{ij}:\cD_j\to\bR$, where $\bR_+=(0,\infty)$.
We use a directed graph in order to denote for each edge a positive flow direction, i.e., if $\{i,j\}\in\cE$ then $\{j,i\}\not\in\cE$, which leads to the identity $\vphi_{ij}(t,x_{ij})=-\vphi_{ji}(t,L_{ij}-x_{ij})$.
We define the set of controllers $\cC\subset \cE\times\{+,-\}$, where $\{i,j\}\equiv\{i,j,+\}\in\cC$ denotes a controller located at node $i\in\cV$ that augments the density of gas flowing into edge $\{i,j\}\in\cE$ in the $+$ direction, while $\{j,i\}\equiv\{i,j,-\}\in\cC$ denotes a controller located at node $j\in\cV$ that augments density into edge $\{i,j\}\in\cE$ in the negative direction.  Compression is then modeled as a multiplicative ratio $\alpha_{ij}:\cT\to\bR_+$ for $\{i,j\}\in\cC$.

We denote by $s_j:\cT\to\bR$ the density of gas entering the network from a node $j\in\cV_S$, where the set $\cV_S$ denotes large supply terminals we call ``slack'' junctions, able to supply any mass flux at the given density.  A mass flux withdrawal (or injection, if negative) at a junction $j\in\cV_D=\cV\setminus V_S$ is denoted by $d_j:\cT\to\bR$, where $\cV_D$ is the set of demand (non-``slack'') nodes.  The functions $\{\alpha_{ij}\}_{\{i,j\}\in\cC}$, $\{d_j\}_{j\in\cV_D}$, and $\{s_j\}_{j\in\cV_S}$ create nodal balance conditions of the form
\begin{align}
\alpha_{ji}(t)\rho_{jk}(t,0) =\alpha_{jk}(t)\rho_{ij}(t,L_j), \qquad \qquad \qquad \qquad \nonumber \\
\qquad \qquad \qquad \qquad \fA j\in\cV_D \text{ and } \{i,j\},\{j,k\}\in\cE,  \label{bc1a}\\
d_j(t) =\sum_{i\in\cV_D}\vphi_{ij}(t,L_{ij})- \sum_{k\in\cV_D}\vphi_{jk}(t,0), \,\,\, \fA j\in\cV_D,  \label{bc1b} \\
\rho_{ij}(t,0)=s_i(t), \quad \fA i\in\cV_S \label{bc1c}
\end{align}

To more conveniently represent the dynamics \eqref{euler1a}-\eqref{euler1b} for each edge, we apply the dimensional transformations
\begin{align} \label{nondim}
\hat{t}=\frac{t}{\ell/a}, \quad \hat{x}=\frac{x}{\ell}, \quad \hat{p}=\frac{\rho}{\rho_0}, \quad \hat{\vphi}=\frac{\vphi}{a\rho_0},
\end{align}
where $\ell$ and $\rho_0$ are nominal length and density, to yield the non-dimensional edge dynamics
\begin{align}
 \dS \pp_t\rho_{ij}+\pp_x\vphi_{ij} & = 0, & \fA \{i,j\}\in\cE \label{euler2a} \\
\dS\pp_t\vphi_{ij} + \pp_x \rho_{ij} & = -\frac{\lambda_{ij}\ell}{2D_{ij}}\frac{\vphi_{ij} |\vphi_{ij}|}{\rho_{ij}}, & \fA \{i,j\}\in\cE \label{euler2b}
\end{align}
in which the hats have been omitted for readability.  We henceforth use the non-dimensional units, and consider \eqref{euler2a}-\eqref{euler2b} to represent flow dynamics on a pipeline segment.  Observe that this non-dimensionalization is not edge-dependent, hence the same factors $\ell$ and $\rho_0$ are used system-wide.

With this notation, we formulate two PDE-constrained OCPs for gas pipeline networks, for which the edge dynamics  \eqref{euler2a}-\eqref{euler2b} and nodal conditions \eqref{bc1a}-\eqref{bc1c} form dynamic constraints.  Each is subject to transient withdrawals $d_j(t)$ for $j\in\cV_D$, available supply densities $s_j(t)$ for $j\in\cV_S$, and box constraints on the density and compression of the form
\begin{align}
\rho_{ij}^{\min}\leq \rho_{ij}(t,x_{ij}) \leq \rho_{ij}^{\max}, \quad \fA \{i,j\}\in\cE,  \label{boxcon1a}\\
1\leq \alpha_{ij}(t) \leq \alpha_{ij}^{\max}, \quad \fA \{i,j\}\in\cC. \label{boxcon1b}
\end{align}
For simplicity, we choose terminal conditions on the state and control variables to be time-periodic,
\begin{align}
\rho_{ij}(0,x_{ij})&=\rho_{ij}(T,x_{ij}),  &\fA \{i,j\}\in\cE, \label{termcon1a}\\
\phi_{ij}(0,x_{ij})&=\phi_{ij}(T,x_{ij}),  &\fA \{i,j\}\in\cE, \label{termcon1b}\\
\alpha_{ij}(0)&=\alpha_{ij}(T), &\fA \{i,j\}\in\cC, \label{termcon1c}
\end{align}
though if the initial conditions are specified we may use
\begin{align}
\sum_{\{i,j\}\in\cE}\int_0^{L_{ij}}\rho_{ij}(0,x)\rd x = \sum_{\{i,j\}\in\cE}\int_0^{L_{ij}}\rho_{ij}(T,x)\rd x, \label{termcon2}
\end{align}
which is a periodic condition on total mass in the system.

The ETC objective function, which aggregates compression costs of the form \eqref{obj1} throughout the network, is
\begin{align}
\!\! J_E=\!\!\!\sum_{\!\!\!\!\{i,j\}\in\cC}\! \int_0^T \!\!\frac{|\phi_{ij}(t,0)|}{\eta_{ij}}\bp{(\max\{\alpha_{ij}(t),1\})^{2m}-1}  \rd t, \label{objetc}
\end{align}
where $\eta_{ij}$ is the efficiency of compressor $\{i,j\}\in\cC$.  We formulate the ETC OCP for day-ahead operational planning, in which the system state is expected to repeat on a 24-hour period $\cT$.  The decision functions are nodal controls $\alpha_{ij}(t)$ for $\{i,j\}\in\cC$, and the problem is
\begin{equation} \label{prob:etc}
\begin{array}{llll}
\min & J_E \text{ in } \eqref{objetc} \\ s.t.
& \text{dynamic constraints: } \eqref{euler2a}-\eqref{euler2b}  \\
& \text{nodal conditions: } \eqref{bc1a}-\eqref{bc1c}\\
& \text{density \& control constraints: } \eqref{boxcon1a}-\eqref{boxcon1b} \\
& \text{terminal constraints: } \eqref{termcon1a}-\eqref{termcon1c}
\end{array}
\end{equation}

In the case that the ETC OCP \eqref{prob:etc} does not have a feasible solution, we wish to determine the control protocol that comes closest to fulfilling the desired gas deliveries $d_j(t)$.  In current practice, natural gas is purchased in day-ahead contracts for firm or non-firm delivery.  Customers of the former type are guaranteed deliveries, while the latter may be cut off if their demand profile is thought to impact pipeline security,  i.e. the density constraints \eqref{boxcon1a}.  We therefore formulate the MLS objective function in order to minimize the gap between the desired and achievable deliveries to non-firm customers at the nodes $\cM=\{j_1,\ldots,j_{\Sigma}\}\subset \cV_D$.  Suppose that for $j\in\cM$, the desired withdrawal from the network is $d_j^*(t)$, and the MLS objective is given by
\begin{align}
J_{M}=\sum_{j\in\cM} \int_0^T c_j(t)(d_j(t)-d_j^*(t))^2 \rd t. \label{objmls}
\end{align}
Here, $J_M$ is the sum of $L_2$ norms of delivery gaps weighted by the marginal quadratic costs $c_j(t)$ of load-shedding at nodes $j\in\cM$. Minimizing $J_M$ leads to the desired distribution of load-shedding across the pipeline system.

The decision functions for the MLS OCP are nodal controls $\alpha_{ij}(t)$ for $\{i,j\}\in\cC$, as well as positive withdrawals $d_j(t)$ for $j\in\cM$, for which we must provide constraints
\begin{align} \label{delivconst1}
d_j(0)=d_j(T) \text{ and } d_j(t) \geq 0, \,\, \fA j\in\cM.
\end{align}
The MLS problem is given by
\begin{equation} \label{prob:mls}
\begin{array}{llll}
\min & J_M \text{ in } \eqref{objmls} \\ s.t.
& \text{dynamic constraints: } \eqref{euler2a}-\eqref{euler2b}  \\
& \text{nodal conditions: } \eqref{bc1a}-\eqref{bc1c}\\
& \text{density \& control constraints: } \eqref{boxcon1a}-\eqref{boxcon1b} \\
& \text{terminal constraints: } \eqref{termcon1a}-\eqref{termcon1c} \\
& \text{delivery constraints: } \eqref{delivconst1}
\end{array}
\end{equation}
The OCPs \eqref{prob:etc} and \eqref{prob:mls} are in general analytically intractable, and pose challenges for computational solution even for $|\cE|=1$, i.e., a pipeline with no controllers or withdrawals except at the boundaries \cite{herty14}.  We proceed to derive the RNF control system based on a model reduction technique for flows on gas networks \cite{grundel13a}.  The RNF replaces the dynamic constraints in the OCPs \eqref{prob:etc} and \eqref{prob:mls}.


\section{Control System Model Reduction Of Gas Network Dynamics} \label{sec:model}

 \begin{figure}[t]
\includegraphics[width=\linewidth]{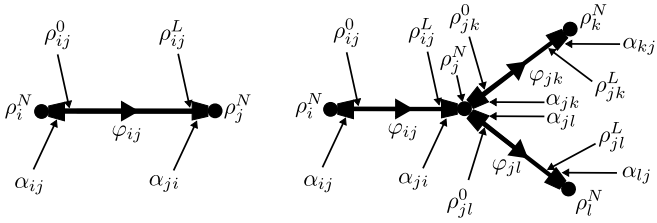} \caption{Relation \eqref{vertmap} between nodal densities $\rho_j^N$ and endpoint densities $\rho_{ij}^0$ and $\rho_{ij}^L$, illustrated for a single edge (left) and for a joint (right).} \vspace{1ex} \label{fig:net1}
\end{figure}

We construct a control system using a tractable yet accurate model of the network flow dynamics \eqref{euler2a}-\eqref{euler2b} with nodal conditions \eqref{bc1a}-\eqref{bc1c}, extending previous modeling work \cite{grundel13a}.  A pipeline segment of (non-dimensional) length $L$ is modeled as a lumped element by integrating \eqref{euler2a}-\eqref{euler2b} with respect to $x$,
\begin{align}
\dS\int_0^L(\pp_t\rho+\pp_x\vphi)\rd x&=0, \label{eulerint1} \\
\dS\int_0^L(\pp_t\vphi + \pp_x\rho)\rd x &= -\frac{\lambda\ell}{2D}\int_0^L \frac{\vphi |\vphi|}{\rho}\rd x, \label{eulerint2}
\end{align}
and then evaluating these integrals of $\pp_t$, $\pp_x$, and nonlinear terms using the trapezoid rule, the fundamental theorem of calculus, and averaging variables, respectively.  This yields
\begin{align}
\!\!\frac{L}{2}(\dot{\rho}^0+\dot{\rho}^L) &= \vphi^0-\vphi^L, \label{eulerdisc1}\\
\!\!\frac{L}{2}(\dot{\vphi}^0+\dot{\vphi}^L) &= \rho^0-\rho^L- \frac{\lambda\ell L}{4D} \frac{(\vphi^0+\vphi^L) |\vphi^0+\vphi^L|}{\rho^0+\rho^L}, \label{eulerdisc2}
\end{align}
where $\rho^0$, $\vphi^0$ and $\rho^L$, $\vphi^L$ denote density and flux at the start and end of an edge, respectively, and $\dot{\rho}=\ddx{t}\rho$.

We extend \eqref{eulerdisc1}-\eqref{eulerdisc2} to a network by defining input and output flows $\vphi_{ij}^0$, $\vphi_{ij}^L$ and densities $\rho_{ij}^0>0$, $\rho_{ij}^L>0$ for each edge $\{i,j\}\in\cE$.  Then equations \eqref{eulerdisc1}-\eqref{eulerdisc2} and nodal conditions \eqref{bc1a}-\eqref{bc1c} then reduce to a differential-algebraic equation (DAE) system where the edge dynamics are
\begin{align}
\!\!\!\!\! \frac{\dot{\rho}_{ij}^L+\dot{\rho}_{ij}^0}{2} & = \dS -\frac{\vphi_{ij}^L-\vphi_{ij}^0}{L_{ij}}, \label{dae0c}
\\ \!\!\!\!\!\dS  \frac{\dot{\vphi}_{ij}^L+\dot{\vphi}_{ij}^0}{2} & \dS  = -\frac{\rho_{ij}^L-\rho_{ij}^0}{L_{ij}}- \frac{\lambda_{ij} \ell}{4D_{ij}}\tS\bp{\frac{(\vphi_{ij}^L+\vphi_{ij}^0) |\vphi_{ij}^L+\vphi_{ij}^0|}{\rho_{ij}^L+\rho_{ij}^0} }, \!\label{dae0d}
\end{align}
and the nodal constraints become
\begin{align}
\!\! 0&=\alpha_{jk}\rho_{ij}^L-\alpha_{ji}\rho_{jk}^0, \, \left\{\begin{array}{l}\fA j\in\cV_D, \\  \{i,j\},\{j,k\}\in\cE, \end{array}\right. \label{dae0a}\\
d_j& = \dS \sum_{i\in\cV_D}\vphi_{ij}^L - \sum_{k\in\cV_D}\vphi_{ij}^0, \,\,\, \fA j\in\cV_D,  \label{dae0b}\\
\rho_{ij}^0& =s_i, \fA i\in\cV_S. \label{dae0e}
\end{align}
Here \eqref{dae0a} represents continuity of pressure at junctions with jumps in the case of compression or regulation, \eqref{dae0b} represents flow balance at junctions, and \eqref{dae0c}-\eqref{dae0d} represent flow dynamics on each segment.

We apply the graph notation in Section \ref{sec:netflow} to express \eqref{dae0c}-\eqref{dae0e} in matrix-vector form.  Suppose that $V=|\cV|$, and assign to each node an index in $[V]$, where $[N]=\{1,\ldots,N\}$ for a positive integer $N\in\bN$.  Define mappings $\pi_v^0:\cE\to[V]$ and $\pi_v^L:\cE\to[V]$, which map the first and last vertices of an edge to the vertex ordering.  Also suppose that $E=|\cE|$ and assign to each edge an index in $[E]$, and define  $\pi_e:\cE\to[E]$, which maps the edges to this ordering.
We now assign to each node in the range of $\pi_v$ a unique internal density, and write a total nodal density state vector $\rho^N=(\rho_1^N,\ldots,\rho_V^N)^T$.  The mappings $\pi_v^0$ and $\pi_v^L$ can be used to state the dependence of $\rho^N$ on the variables $\{\rho_{ij}^0\}$ and $\{\rho_{ij}^L\}$ and control variables,
\begin{align} \label{vertmap}
\rho_{ij}^0=\alpha_{ij}\rho_{\pi_v^0(ij)}^N \text{  and  } \rho_{ij}^L=\alpha_{ji}\rho_{\pi_v^L(ij)}^N.
\end{align}
Equation \eqref{vertmap} will be used to state \eqref{dae0c}-\eqref{dae0d} in terms of nodal densities $\rho^N$.  We also write state vectors $\vphi^0=(\vphi_1^0,\ldots,\vphi_E^0)^T$ and $\vphi^L=(\vphi_1^L,\ldots,\vphi_E^L)^T$, where $\vphi_k^0$ and $\vphi_k^L$ are indexed by $k=\pi_e(ij)$.

We then define the time-dependent weighted incidence matrix $B:\bR^{E}\to\bR^{V}$ by
\begin{align} \label{incidence0}
B_{ik} = \left\{ \begin{array}{ll}  \alpha_{ij} & \text{edge $k=\pi_e(ij)$ enters node $i$,} \\ -\alpha_{ij} & \text{edge $k=\pi_e(ij)$ leaves node $i$,} \\ 0 & \text{else} \end{array}\right.
\end{align}
as well as the incidence matrix $A=\rM{sign}(B)$.  We define the collection of withdrawal fluxes $d=(d_1,\ldots,d_M)^T$ with $M=|\cV_D|$, where $d_k$ is negative if an injection. Also define the slack node densities as $s=(s_1,\ldots,s_b)^T=\{\rho^N_j\}_{j\in\cV_S}$, where $b=|\cV_S|$, and demand node densities as $\rho=(\rho_1,\ldots,\rho_M)^T=\{\rho^N_j\}_{j\in\cV_D}$, so that $b+M=V$.  Then let $A_s,B_s\in\bR^{b\times E}$ denote the submatrices of rows of $A$ and $B$ corresponding to $\cV_S$, and let $A_d,B_d\in\bR^{M\times E}$ correspond similarly to $\cV_D$.  Next, let $A_L$ and $A_0$ denote the positive and negative parts of $A_d$, so that $A_d=A_L+A_0$.  Also, define the diagonal matrices $\Lambda,K\in\bR^{E\times E}$ by $\Lambda_{kk}=L_k$ and $K_{kk}=\ell\lambda_k/D_k$, where $L_k$, $\lambda_k$, and $D_k$ are the nondimensional length, friction coefficient, and diameter of edge $k=\pi_e(ij)$.       Finally, define a function $g:\bR^E\times\bR_+^E\to\bR^E$ by $g_j(x,y)=x_j|x_j|/y_j$.  Then \eqref{dae0c}-\eqref{dae0e} can be written
\begin{align}
\!\!\!d& = \dS A_L\vphi_L + A_0\vphi_0,  \label{dae1a}\\
\!\!\!\dS  |B_s^T|\dot{s} + |B_d^T|\dot{\rho} & = \dS -4\Lambda^{-1} \thalf(\vphi_L-\vphi_0), \label{dae1b}
\\ \!\!\!\dS  \thalf(\dot{\vphi}_L+\dot{\vphi}_0) &= -\Lambda^{-1}(B_s^Ts+B_d^T\rho) \nonumber\\ & \quad -Kg(\thalf(\vphi_L+\vphi_0),|B_s^T|s+|B_d^T|\rho) \label{dae1c}
\end{align}
Furthermore, note that $A_d=A_L+A_0$ and $|A_d|=A_L-A_0$, so by defining $\vphi=\thalf(\vphi_L+\vphi_0)$ and $\vphi_-=\thalf(\vphi_L-\vphi_0)$ we may replace \eqref{dae1a} with $d = A_d\vphi + |A_d|\vphi_-$, the right hand side of \eqref{dae1b} with $-4\Lambda^{-1} \vphi_-$, and the left hand side of \eqref{dae1c} with $\dot{\vphi}$.  Then multiplying \eqref{dae1b} by $|A_d|\Lambda$ results in $|A_d|\Lambda|B_s^T|\dot{s} + |A_d|\Lambda|B_d^T|\dot{\rho} = -4|A_d|\vphi_- = 4(A_d\vphi-d )$.  The DAE system \eqref{dae1a}-\eqref{dae1c} may then be written as an ODE
\begin{align}
\!\!\! \dS  \dot{\rho} & = \dS (|A_d|\Lambda|B_d^T|)^{-1}[4(A_d\vphi - d) - |A_d|\Lambda|B_s^T|\dot{s}], \label{ode0a}
\\ \!\!\! \dS  \dot{\vphi} &= -\Lambda^{-1}(B_s^Ts+B_d^T\rho) -Kg(\vphi,|B_s^T|s+|B_d^T|\rho), \label{ode0b}
\end{align}
where $\rho$ are nodal densities and $\vphi$ approximate mass flux on the edges.  For a connected graph, $A_d\in\bR^{M\times E}$ and $B_d\in\bR^{M\times E}$ are full rank, and therefore $|A_d|\Lambda|B_d^T|$ is invertible.  Time-varying parameters are gas withdrawals $d\in\bR^M$, input densities $s\in\bR_+^b$, and compressions/regulations $\alpha_{ij}\in\cC$.
We now give two consistency results.  The RNF for a pipeline reduces to \eqref{euler2a}-\eqref{euler2b} as the maximum spatial discretization step approaches zero, and reduces to the Weymouth equations \cite{borraz10phd} in the steady-state.

\begin{prop}{The RNF \eqref{ode0a}-\eqref{ode0b} is a consistent spatial discretization of the PDE \eqref{euler2a}-\eqref{euler2b} for a pipeline of dimensional length L, modeled as a chain of $m$ segments of uniform length $L/m$, and which has no compressors or intermediate withdrawals ($d_j\equiv0$ for $j\neq m$).} \label{prop:consist} \\
\pf
The pipeline is represented as a graph $\cG$ with $m=E$, $V=E+1$, and each node at the points $x_i=iL/m$ for $i=1,\ldots,m-1$ is connected to two edges and for $i=0,m$ is connected to one edge.   It can be shown that $W=|A_d|\Lambda|B_d^T|$ is tridiagonal with $W_{i,i-1}=L/(m\ell)$, $W_{i,i}=2L/(m\ell)$, and $W_{i,i+1}=L/(m\ell)$, and matrix $A_d$ satisfies $(A_d)_{i,i}=1$, $(A_d)_{i,i+1}=-1$.   With $\rho_i(t)=\rho(t,x_i)$ and $\vphi_i(t)=\vphi(t,x_{i}-\thalf \ell)$ at time $t$,  \eqref{ode0a}-\eqref{ode0b} yield
\begin{align}
&\dS  \frac{1}{4}(\dot{\rho}_{i-1}+2\dot{\rho}_i+\dot{\rho}_{i+1})  = -\frac{L}{m\ell}(\vphi_i-\vphi_{i+1}), \label{consist0a}
\\ &\dS  \dot{\vphi}_i = -\frac{L}{m\ell}(\rho_i-\rho_{i-1}) -\frac{\lambda\ell}{2D}\frac{\vphi_i|\vphi_i|}{\thalf(\rho_i+\rho_{i-1})}, \label{consist0b}
\end{align}
with scalar boundary conditions $\rho_1=s$ and $\vphi_m=d_m$.  Taking $m\to\infty$ yields \eqref{euler1a}-\eqref{euler1b} with the transformations \eqref{nondim}.
\end{prop}

\begin{rem}
Proposition \ref{prop:consist} implies that the solution to \eqref{ode0a}-\eqref{ode0b} using an implicit ODE integrator, which adapts the time-step to maintain stability and accuracy, will converge to the solution to \eqref{euler2a}-\eqref{euler2b} as the space discretization $M$ is increased \cite{verwer84}.  Indeed, a simulation of the RNF approximation \eqref{ode0a}-\eqref{ode0b} for a standard pipeline model \cite{herty10} with slow transients in input pressure and output flux yielded similar results as a solution of \eqref{euler2a}-\eqref{euler2b} using an operator split-step method for hyperbolic PDE systems \cite{zlotnik15dscc}.
\end{rem}

\begin{prop}{When $\dot{\rho}=0$, $\dot{\vphi}=0$, $\dot{s}=0$, $\dot{d}=0$, and $\dot{\alpha}_{ij}=0$ for all $\{i,j\}\in\cC$, equations \eqref{ode0a}-\eqref{ode0b} reduce to the steady-state balance laws for a gas network \cite{misra15}.} \label{prop:static} \\
\pf Equation \eqref{ode0a} is reduced to $A_d\vphi=d$, which is nodal conservation of flow on demand nodes.  Recall that $\rho^N$ contains all nodal densities, so equation \eqref{ode0b} leads to
\begin{align}\label{weymouth1}
\Lambda K\vphi\odot|\vphi|&=(B^T\rho^N)\odot (|B^T|\rho^N),
\end{align}
where $\odot$ denotes the point-wise vector product.  Returning to $\{\rho_{ij}^0\}$ and $\{\rho_{ij}^L\}$ from the nodal densities using \eqref{vertmap} and reverting to the dimensional variables leads to the Weymouth equations for static gas networks \cite{borraz10phd},
\begin{align}\label{weymouth2}
-\frac{\lambda_{ij} ZRT L_{ij}}{D_{ij}}\vphi_{ij}|\vphi_{ij}|&= (\rho_{ij}^L)^2- (\rho_{ij}^0)^2.
\end{align}
\end{prop}
\begin{rem}
The potential equation \eqref{weymouth2}, where $\{\rho_{ij}^0\}$ for any slack node $i\in\cV_S$ is given, together with the relation $A_d\vphi=d$, characterize the steady-state balance laws \cite{misra15}.
\end{rem}
With the above results, we see that the RNF \eqref{ode0a}-\eqref{ode0b} can be used to represent both the PDE model \eqref{euler1a}-\eqref{euler1b} for a dynamic network and the static equilibrium equations.

To implement the RNF in the transient ETC \eqref{prob:etc} and MLS \eqref{prob:mls} OCPs, we rewrite the constraints and objective functions in Section \ref{sec:netflow} in terms of 
the nodal densities $\rho$ and edge fluxes $\vphi$ used in \eqref{ode0a}-\eqref{ode0b}.  The inequality constraint \eqref{boxcon1a} becomes
\begin{align}
&\rho_i^{\min}\leq \alpha_{ij}(t)\rho_i(t) \leq \rho_i^{\max},  \label{boxcon2a}
\end{align}
and the time-periodic terminal conditions \eqref{termcon1a}-\eqref{termcon1b} become
\begin{align}
&\rho(0)=\rho(T), \,\,\,\, \vphi(0)=\vphi(T),  \label{termcon3a}
\end{align}
and total mass conservation \eqref{termcon2} becomes
\begin{align}
\!\!\!\! \Bone^T\Lambda(|B_s^T|(s(0)-s(T))+|B_d^T|(\rho(0)-\rho(T)))=  0 \label{termcon4}
\end{align}
where $\Bone$ is a vector with one in each coordinate. The ETC objective function \eqref{objetc} becomes
\begin{align}
\!\!\!\!\! J_E=\!\!\!\sum_{\!\!\!\!\{i,j\}\in\cC}\! \int_0^T \!\!\frac{|\vphi_{\pi_e(ij)}(t)|}{\eta_{ij}}\bp{(\max\{\alpha_{ij}(t),1\})^{2m}\!-\!1}  \rd t \label{objetc2}
\end{align}
while the objective \eqref{objmls} for MLS remains unchanged.  Using the RNF, the reduced ETC OCP takes the form
\begin{equation} \label{prob:etc-rm}
\begin{array}{llll}
\min & J_E \text{ in } \eqref{objetc2} \\ s.t.
& \text{RNF constraints: } \eqref{ode0a}-\eqref{ode0b}  \\
& \text{density \& control constraints: } \eqref{boxcon2a}, \eqref{boxcon1b} \\
& \text{terminal constraints: } \eqref{termcon3a}, \eqref{termcon1c}
\end{array}
\end{equation}
The reduced MLS OCP is given as
\begin{equation} \label{prob:mls-rm}
\begin{array}{llll}
\min & J_M \text{ in } \eqref{objmls} \\ s.t.
& \text{RNF constraints: } \eqref{ode0a}-\eqref{ode0b}  \\
& \text{density \& control constraints: } \eqref{boxcon2a}, \eqref{boxcon1b} \\
& \text{terminal constraints: } \eqref{termcon3a}, \eqref{termcon1c} \\
& \text{delivery constraints: } \eqref{delivconst1}
\end{array}
\end{equation}
Although we have reduced the instantaneous states to the vectors $\rho(t)$ and $\vphi(t)$ from the continuous functions $\rho_{ij}(t,x_{ij})$ and $\vphi_{ij}(t,x_{ij})$ used in Section \ref{sec:netflow}, the OCPs \eqref{prob:etc-rm} and \eqref{prob:mls-rm} require optimization on the space of functions $\alpha_{ij}(t)$ for all $\{i,j\}\in\cC$.  Therefore, we employ a method for time-discretization of such problems to NLPs.

\section{Pseudospectral Optimal Control} \label{sec:pseudospectral}

\noindent We review a numerical scheme for transcribing OCPs into NLPs using pseudospectral discretization.  Consider the OCP
\begin{align}
	\min_u\ \ & J(x,u)=\int_0^T \mathcal{L}(t,x(t),u(t))dt, \label{eq:ocp0a} \\
	{\rm s.t.}\ \ & \dot{x}(t) = f(t,x(t),u(t)), \label{eq:ocp0b} \\
	& e(x(0),x(T)) = 0, \label{eq:ocp0c} \\
	& g(x(t),u(t))\leq 0, \label{eq:ocp0d}
\end{align}
on $\cT=[0,T]$.  Here $\mathcal{L}\in C^\kappa$ is in the space $C^\kappa$ of continuous functions with $\kappa$ classical derivatives, and the dynamic constraints $f\in C_n^{\kappa-1}$ are in the space $C_n^{\kappa-1}$ of $n$-vector valued $C^{\kappa-1}$ functions, with respect to the state, $x(t)\in\bR^n$, and control, $u(t)\in\bR^m$.  The functions $e$ and $g$ are terminal and path constraints, respectively.   Finally, the admissible set for controls $u$ includes the piecewise $C_m^{\kappa}$ functions on $\cT$.  For details, refer to \cite{ruths11cdc}.  We now derive a direct collocation procedure for constructing a finite-dimensional nonlinear program that approximates the problem \eqref{eq:ocp0a}-\eqref{eq:ocp0d}.

We use Lagrange polynomial interpolation to approximate the states $x$ and controls $u$ at a set of collocation points $t_k$:
\begin{eqnarray}
	\label{eq:Ix} &x(t) \approx  \wh{x}_N(t) = \sum_{k=0}^N \bar{x}_k \ell_k(t), \\
	\label{eq:Iu} &u(t) \approx \wh{u}_N(t) = \sum_{k=0}^N \bar{u}_k \ell_k(t).
\end{eqnarray}
Lagrange interpolation polynomials satisfy $\ell_k(t_i) = \delta_{ki}$, where $\delta_{ki}$ is the Kronecker delta function \cite{canuto06}.  It follows that $x(t_k) = \wh{x}_N(t_k) = \bar{x}_k$ and $u(t_k) = \wh{u}_N(t_k) = \bar{u}_k$, so the physical meaning of the interpolating polynomial coefficients $\bar{x}_k$ and $\bar{u}_k$ are the values of the state and control variables at the collocation points. Those points are chosen so that the integral in \eqref{eq:ocp0a} and the derivative in \eqref{eq:ocp0b} are computed accurately.  The former is approximated using Legendre-Gauss quadrature, leading to the choice of Legendre polynomials as the orthogonal basis.  Furthermore, the Legendre-Gauss-Lobatto (LGL) quadrature points are chosen to include endpoints of the interval, permitting the terminal constraints to be specified within the scheme.  The LGL quadrature rule for a function $f:[-1,1]\to\bR$ is
\begin{equation}\label{eq:gaussquad}
	\int_{-1}^{1} f(t) dt \approx \sum_{i=1}^{N} f(t_i) w_i, \qquad w_i = \int_{-1}^1 \ell_i(t)dt,
\end{equation}
and is exact if $f\in\mathbb{P}_{2N-1}$ and the nodes $t_i\in\Gamma^{LGL}$, where $\mathbb{P}_{2N-1}$ denotes the set of polynomials of degree at most $2N-1$ and where $\Gamma^{LGL}=\{t_i:\dot{L}_N(t_i)=0, i=1,\ldots N-1\} \bigcup \{-1,1\}$ are the $N+1$ LGL nodes determined by the derivative of the $N^{\text{th}}$ order Legendre polynomial, $\dot{L}_N(t)$, and the interval endpoints \cite{canuto06}.  Because the OCP is defined on $\cT=[0,T]$, whereas Legendre polynomials form a basis on $[-1,1]$, we re-scale time by $\tilde{t}=(2t-T)/T$.

We re-write the Lagrange interpolating polynomials on the LGL collocation nodes in terms of the Legendre polynomial basis, to provide the scheme with derivative and spectral accuracy properties of orthogonal polynomials. Given $t_k\in\Gamma^{LGL}$, we can express the Lagrange polynomials as \cite{boyd01}
\begin{equation}\label{eq:lagleg}
	\ell_k(t) = \displaystyle\frac{1}{N(N+1)L_N(t_k)}\frac{(t^2-1) \dot{L}_N(t)}{t-t_k}.
\end{equation}
\noindent The derivative of (\ref{eq:Ix}) at $t_j \in \Gamma^{LGL}$ is then
\begin{equation}\label{eq:dinterpsc_j}
	\frac{d}{dt} \wh{x}_N(t_j) = \sum_{k=0}^N \bar{x}_k \dot{\ell}_k(t_j)=\sum_{k=0}^N D_{jk}\bar{x}_k,
\end{equation}
\noindent where $D$ is the constant differentiation matrix with elements $D_{ik}=\dot{\ell}_k(t_i)$.
Using \eqref{eq:Ix}, \eqref{eq:Iu}, \eqref{eq:gaussquad}, and \eqref{eq:dinterpsc_j}, the OCP \eqref{eq:ocp0a}-\eqref{eq:ocp0d} is transcribed as the following NLP, in which the decision variables are the polynomial interpolation coefficient vectors $\bar{x}=(\bar{x}_0,\ldots,\bar{x}_N)$ and $\bar{u}=(\bar{u}_0,\ldots,\bar{u}_N)$:
\begin{align}
\!\!\!\!\!	\min\ \ & \bar{J}(\bar{x},\bar{u})=\sum_{k=0}^N \frac{T}{2}\cL(\bar{x}_k,\bar{u}_k)w_k \label{eq:ocp1a}\\
\!\!\!\!\!	{\rm s.t.}\ \ & \dS\sum_{k=0}^ND_{ik}\bar{x}_k=\frac{T}{2}f(t_i,\bar{x}_i,\bar{u}_i) & \!\!\!\!\! \fA\, i=0,1,\ldots,N \label{eq:ocp1b}\\
\!\!\!\!\!	& e(\bar{x}_0,\bar{x}_N) = 0 \label{eq:ocp1c} \\
\!\!\!\!\!	& g(\bar{x}_k,\bar{u}_k) \leq 0 & \!\!\!\!\! \fA \, k=0,1,\ldots,N \label{eq:ocp1d}
	\end{align}
The solutions to \eqref{eq:ocp1a}-\eqref{eq:ocp1d} converge to extrema of \eqref{eq:ocp0a}-\eqref{eq:ocp0d} as $N\to\infty$ at an exponential rate because of the spectral accuracy of polynomial approximations \cite{ruths11cdc}.

\section{Implementation} \label{sec:implement}

The model reduction in Section \ref{sec:model} is in essence a spatial discretization of the the pipeline flow dynamics in Section \ref{sec:pipeline} on a graph $\cG=(\cV,\cE)$, which incorporates boundary conditions at network nodes.  For the model to accurately approximate the true PDE dynamics, this discretization must be sufficiently fine.  Thus to translate from the PDE-constrained problems \eqref{prob:etc} or \eqref{prob:mls} to the reduced OCPs \eqref{prob:etc-rm} and \eqref{prob:mls-rm}, we first create a modified graph $\hat{\cG}=(\hat{\cV},\hat{\cE})$ by adding nodes such that all edges of $\hat{\cE}$ are shorter than a maximum length $\ell$, which can also be used as the non-dimensional constant.  It has been observed \cite{grundel14}, and we have confirmed \cite{zlotnik15dscc}, that $\ell=10$ km is sufficient to adequately represent transients of interest for typical transmission pipelines.  The graph $\hat{\cG}$ is then used to create the RNF \eqref{ode0a}-\eqref{ode0b}, which is used for \eqref{prob:etc-rm} or \eqref{prob:mls-rm}.  The system state variables are then $\rho\in\bR^M$ and $\vphi\in\bR^E$, where $M=|\hat{\cV}_D|$ and $E=|\hat{\cE}|$, so that $x=(\rho,\vphi)$ is the state vector of interest.  The control variables are $\alpha_{ij}$ for $\{i,j\}\in\cC$, where $C=|\cC|$ is the number of compressors/regulators.  In the case of the MLS problem, the withdrawals $d_j$ for $j\in\cM$ are variables as well, and $\Sigma=|\cM|$ are the number of non-firm loads where shedding may occur.  The control vector $u$ contains all $\alpha_{ij}$ and $d_j$ which are to be determined by the system operator.  The OCPs \eqref{prob:etc-rm} and \eqref{prob:mls-rm} can then be expressed in the form of \eqref{eq:ocp0a}-\eqref{eq:ocp0d}, and then approximated as an NLP using the procedure described in Section \ref{sec:pseudospectral}.  In the NLP \eqref{eq:ocp1a}-\eqref{eq:ocp1d}, each element of the vector-valued functions $x$ and $u$ is expressed using $N+1$ Lagrange interpolation coefficients, which leads to $(M+E+C)\times(N+1)$ and $(M+E+C+\Sigma)\times(N+1)$ variables for the ETC and MLS OCPs, respectively.

To guarantee a smooth, physically relevant solution, we add a penalty on the square of the $L_2$ norms of derivatives of the compression ratios to the objective function:
\begin{align} \label{costmod1}
J_S(\alpha)& = \mu\sum_{\{i,j\}\in\cE}||\dot{\alpha}_{ij}||_2^2=\mu\sum_{\{i,j\}\in\cE}\int_0^T(\dot{\alpha}_{ij}(t))^2\rd t \nonumber \\ & \approx \mu\frac{2}{T}\sum_{\{i,j\}\in\cE}\sum_{m=0}^N \bp{\sum_{k=0}^N D_{mk}\bar{\alpha}_{ijk}}^2w_m,
\end{align}
where $\bar{\alpha}_{ij} = (\bar{\alpha}_{ij0},\ldots,\bar{\alpha}_{ijN})^T$ are interpolation coefficients for compression function $\{i,j\}\in\cC$, and $\mu$ is a relative weight, for which an effective empirical value is $N$.  The cost term $J_S$ is eliminated when discussing objective values.

The optimal control scheme is implemented computationally as follows. First, all system parameters including network structure and constraints as well as interpolation coefficients of time-varying withdrawals and injections are used to build MATLAB functions for the objective, constraints, and their gradients with respect to the decision variables.  These are provided, along with random initial conditions that satisfy inequality constraints, to the interior-point solver IPOPT version 3.11.8 running with the sparse linear solver ma57 \cite{biegler09ipopt}.  Convergence of optimization for the case studies below requires only minutes because the gradients are provided to the solver, and the constraint Jacobian provided to ma57 has under 3\% non-zero entries.

 \begin{figure}[t]
\hspace{1.5ex} \includegraphics[width=.96\linewidth, height=5cm]{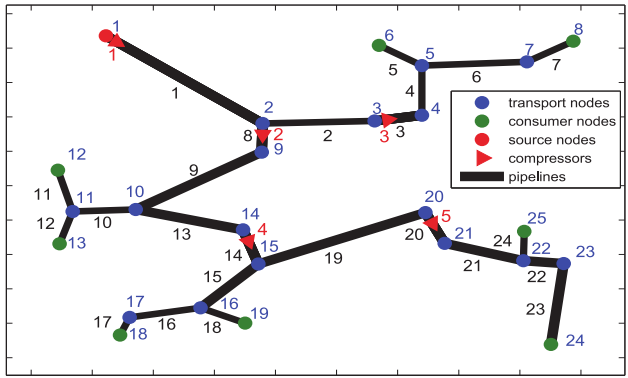} \caption{Example network (not to scale).  Numbering indicates nodes (blue, above/right), edges (black, below/left), and compressors (red, below/right). Thick and thinner lines indicate 36'' and 25'' pipes, respectively.  Pressure is bounded between 500 and 800 psi on all pipes.  Friction factor and sound speed are $\lambda=0.01$ and $a=377.968$ m/s.} \label{fig:net2}
\end{figure}


\section{Examples} \label{sec:examples}

Current industry practice is to assign compression set-points in an ad-hoc manner with the stipulation that gas withdrawals are constant throughout the day \cite{tiso14}.  Because actual gas flows may be up to 80\% above or below the planned-for rates, set-points must be chosen very conservatively, so that pressures may drop far below the rated minimum, which leads to load-shedding and gas price spikes.

Solving the ETC problem addresses such issues by accounting for transient withdrawals, which are usually known by gas system operators on a day-ahead basis \cite{ferc787}.  By utilizing time-dependent dynamical information, the transient ETC is much more likely to have a solution that is actually valid for operations, and load-shedding will be unnecessary.  An intermediate formulation, which we call ``quasi-static'', is examined for comparison, where constant compression set-points are chosen to satisfy pressure constraints given dynamic withdrawals.  When the transient ETC problem does not have a feasible solution, the transient MLS problem can determine the most efficient protocol for load-shedding.  The latter leads to the ultimate utilization of pipeline network capacity while ensuring that pressure limits are not exceeded.

 \begin{figure}[t]
\includegraphics[width=.32\linewidth]{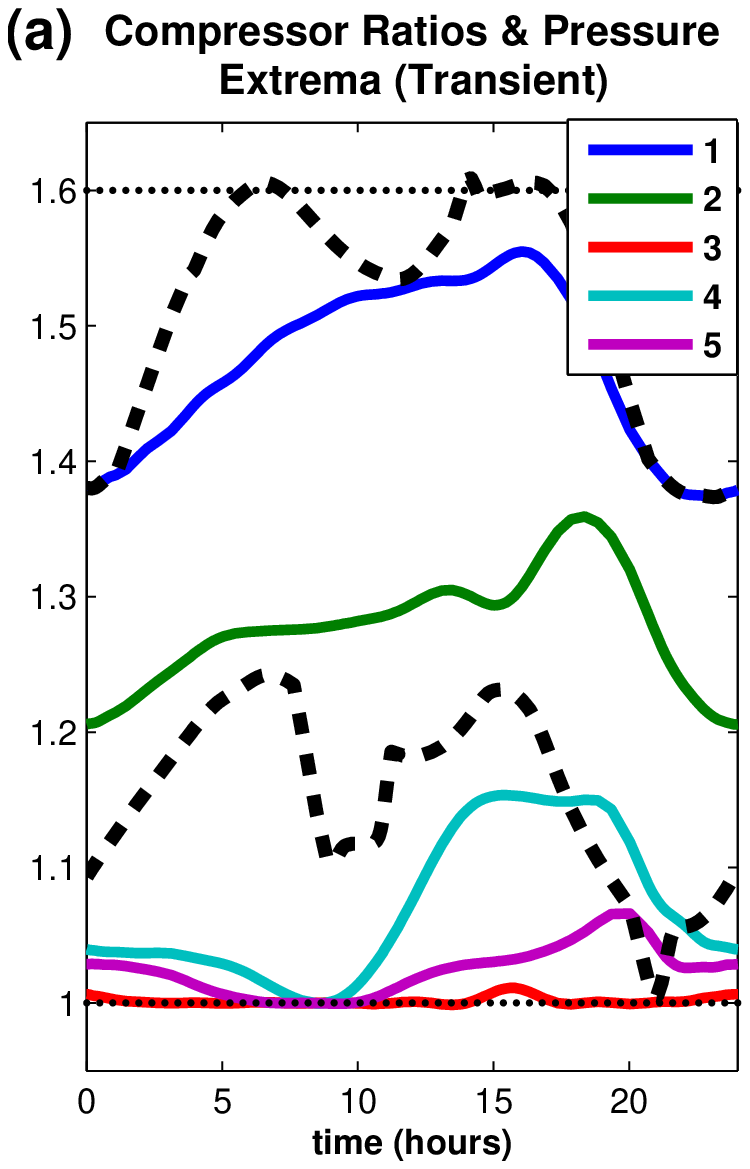} \includegraphics[width=.32\linewidth]{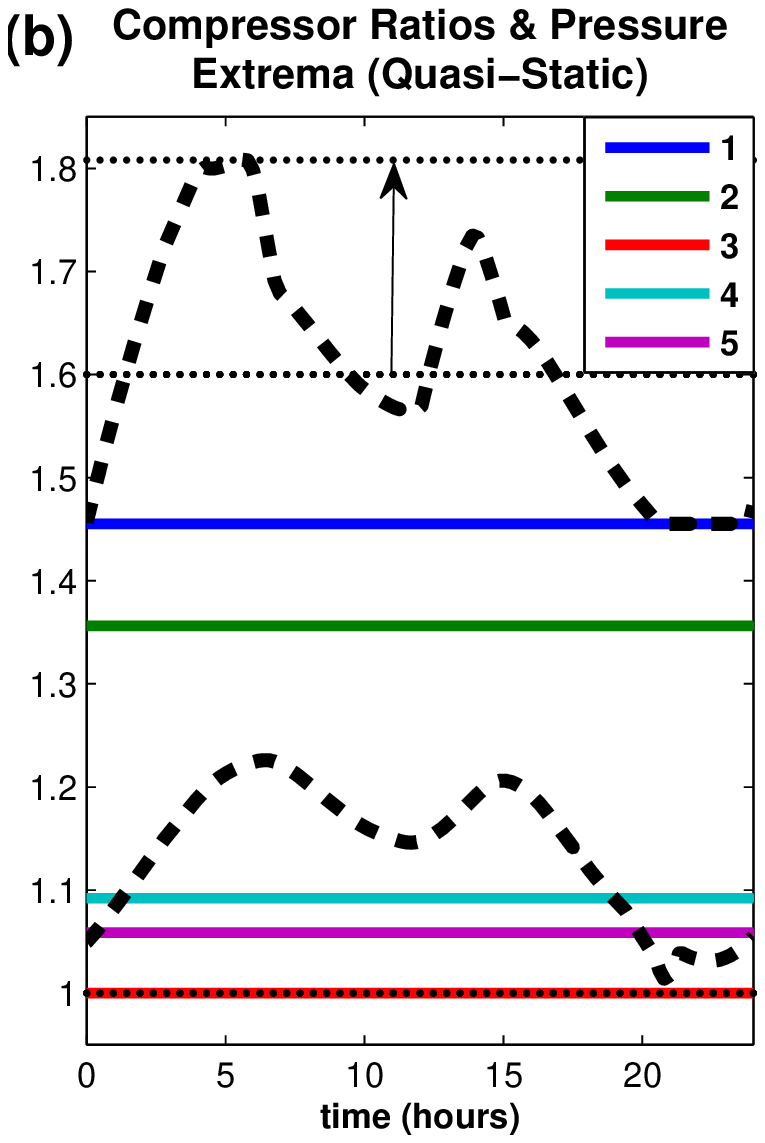}  \includegraphics[width=.32\linewidth]{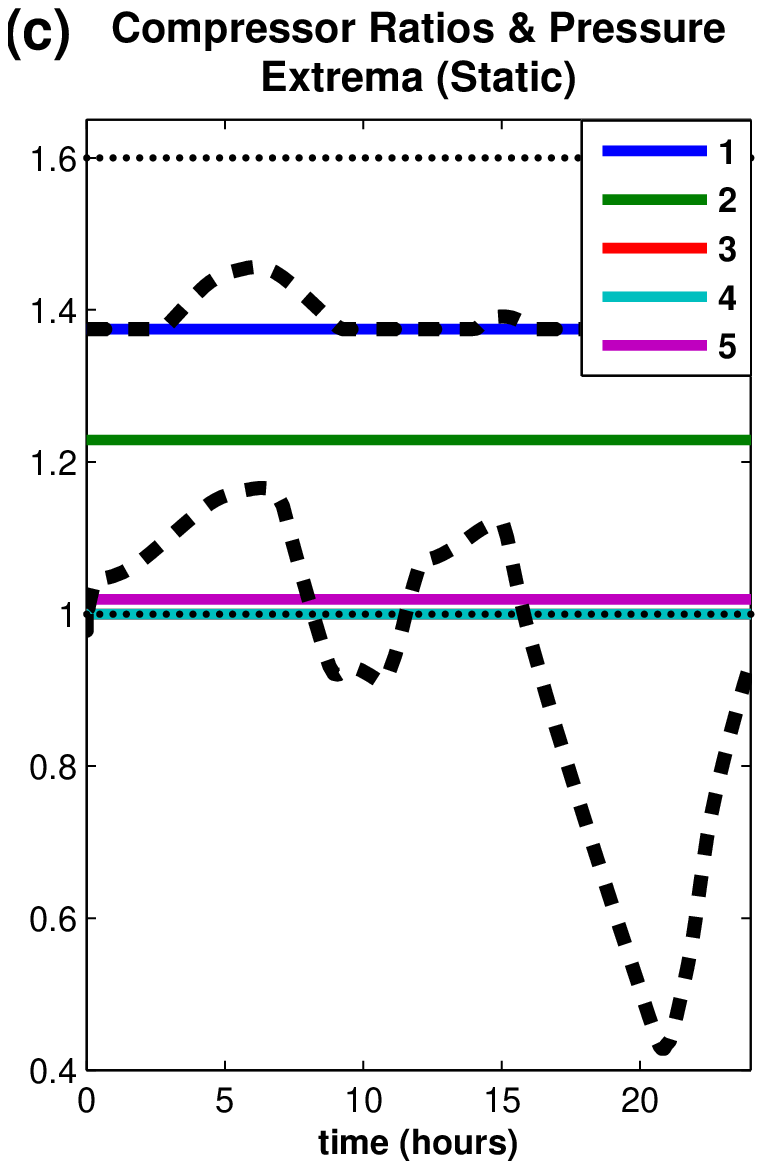} \includegraphics[width=.49\linewidth]{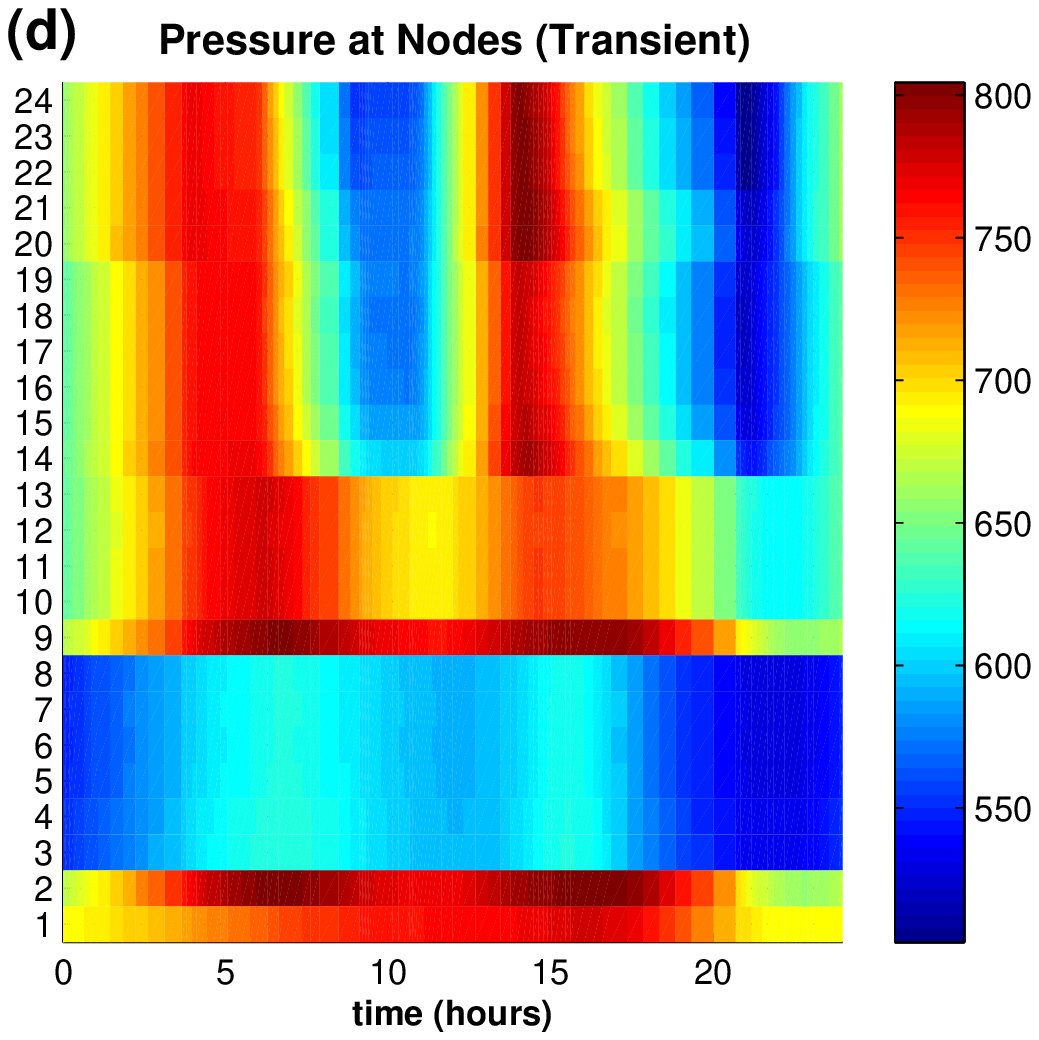} \includegraphics[width=.49\linewidth]{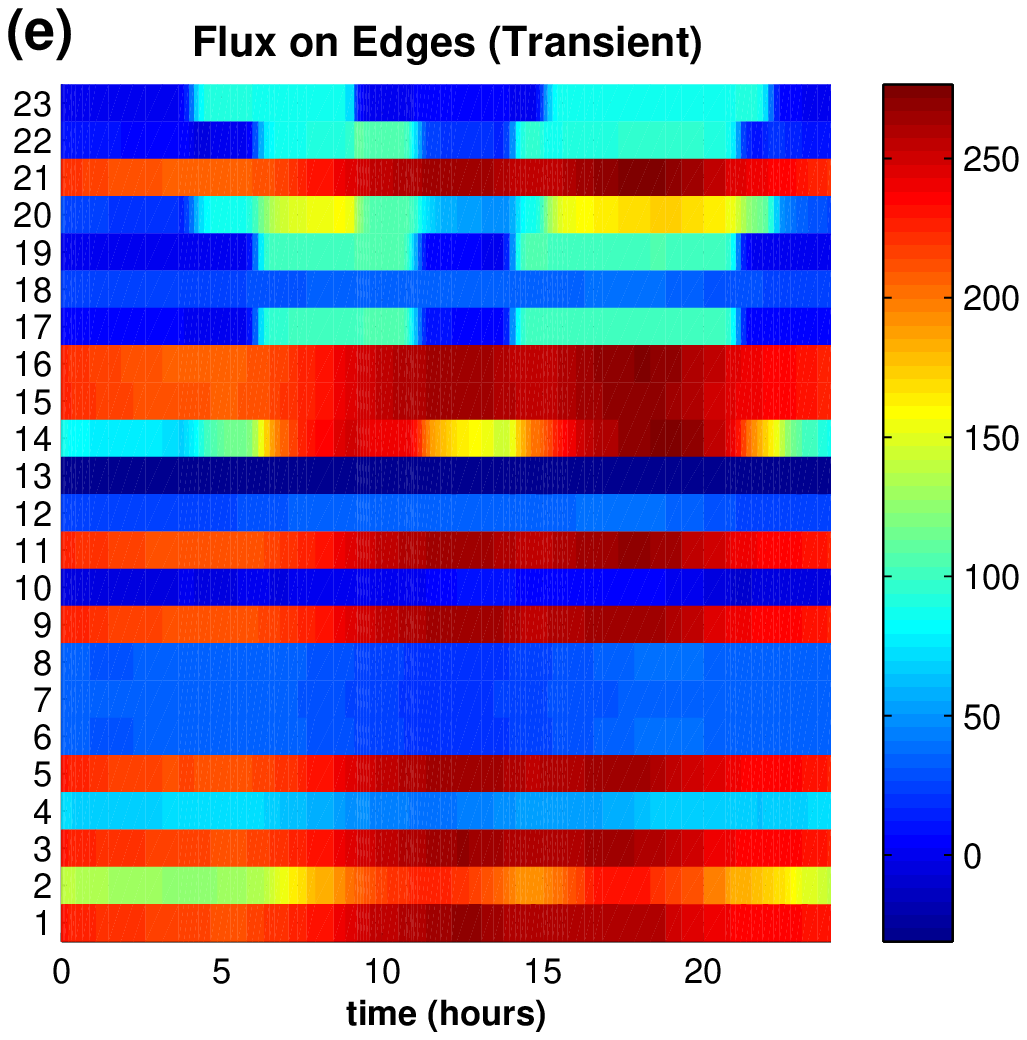}
\includegraphics[width=\linewidth]{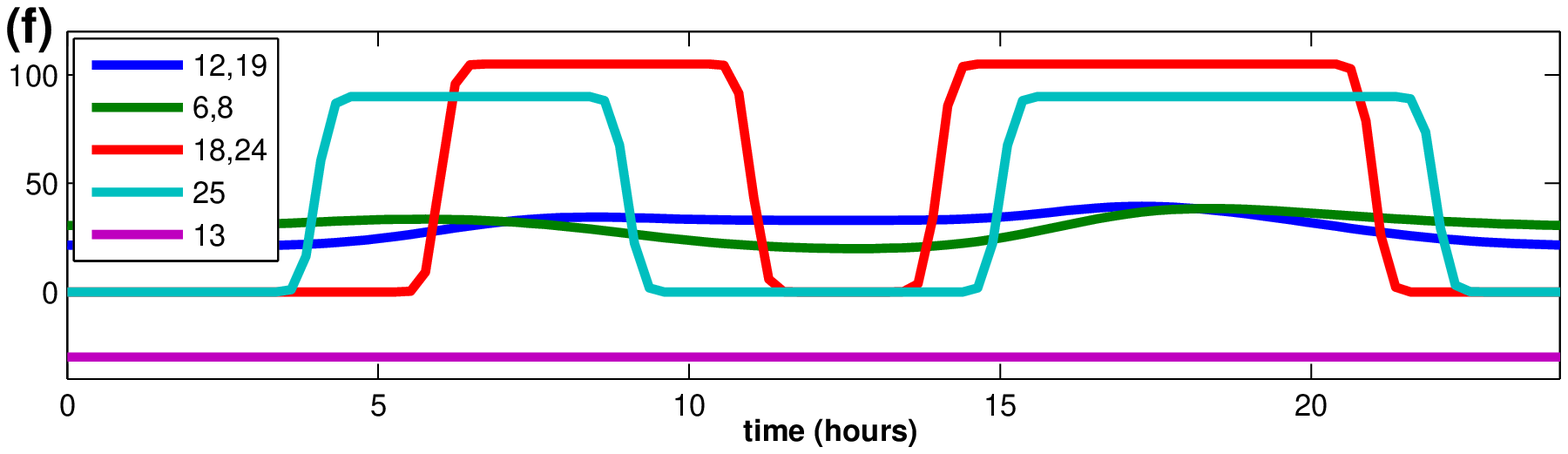}  \caption{ETC Optimization. (a) Transient compression solutions (color), with scaled extremal system pressure (thick dashed line) and limits (dotted line); (b) Quasi-static solution - maximum pressure relaxed by 25\%; (c) Static solution - pressure drop as seen in practice; (d) Pressures at nodes (MPa); (e) Fluxes on edges (kg/m$^2$/s); (f) Withdrawal/injection profiles (kg/m$^2$/s).} \label{fig:net3}
\end{figure}

We consider an example system in Figure \ref{fig:net2} consisting of a tree with $25$ nodes connected by $24$ edges with a total length of 477 km, containing $C=5$ compressors, and a single slack node $j=1$.  For more accurate representation of transients, artificial nodes are added so that pipe segments have maximum length 10km, resulting in $V=62$ nodes (yielding $M=61$ non-slack nodes), and $E=61$ edges.  Gas is supplied at the minimum pressure of 500 psi at the slack node, to be immediately compressed into the network.

{\it ETC Optimization.} The compression ratios for the system in Figure \ref{fig:net2} are optimized to solve the ETC problem \eqref{prob:etc-rm} for the withdrawal profiles in Figure \ref{fig:net3}f.  Two slowly-changing profiles on node sets $\{12,19\}$ and $\{6,8\}$ represent residential and industrial use, while profiles for the sets $\{18,24\}$ and $\{25\}$  represent common single-cycle gas turbine operations.  Node 13 has a constant injection.  Thus the slow and fast demand profiles account for 34.5\% and 65.5\% of the total consumption, respectively.  The problem is solved with $N=25$ time points, so the total number of variables is $(M+E+C)\times(N+1)=3302$, and solution takes under 5 minutes on a laptop computer.  Figure \ref{fig:net3}a shows the transient compression solution, as well as resulting scaled maximum and minimum pressures, which fall within the bounds as indicated.  Figures \ref{fig:net3}b and \ref{fig:net3}c show these results for the quasi-static solution and the fully static solution, respectively.  A feasible quasi-static solution is obtained only when the maximum pressure constraint is relaxed with a 25\% increase. The scaled extremal pressures that result when these solutions are applied to simulations with the transient withdrawals dramatically violate the desired limits, as observed in current pipeline operations \cite{tiso14}.

\begin{figure}[t]
\includegraphics[width=.49\linewidth]{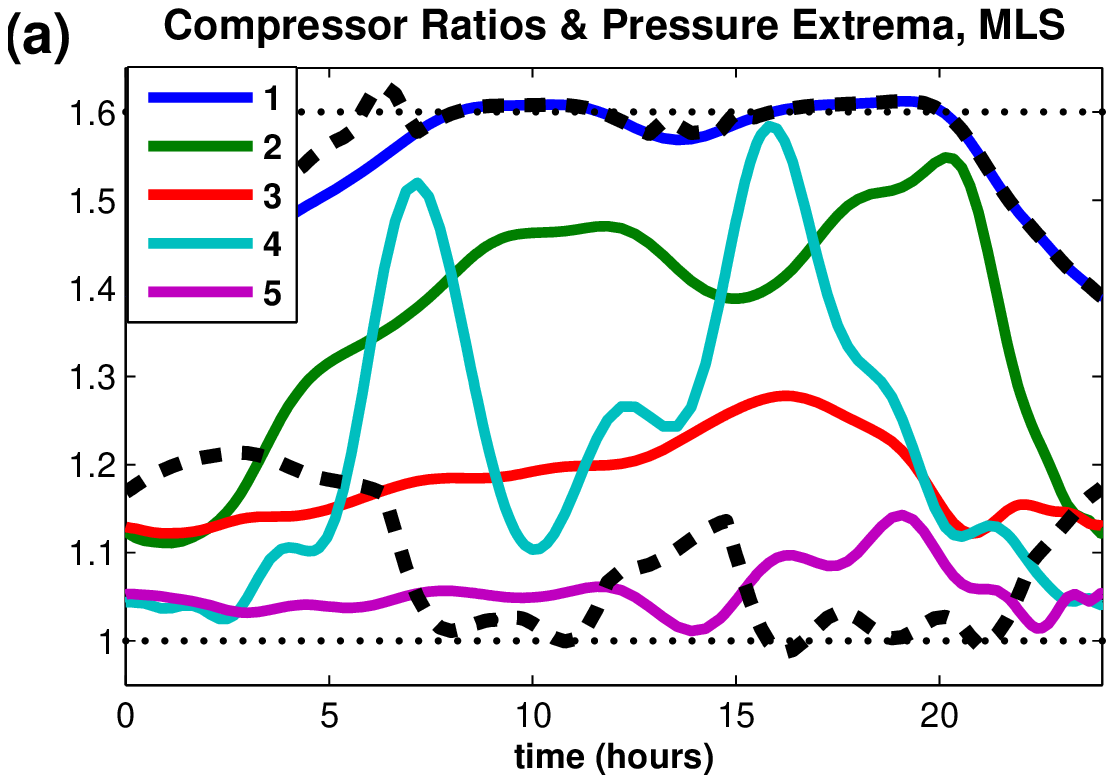} \includegraphics[width=.49\linewidth]{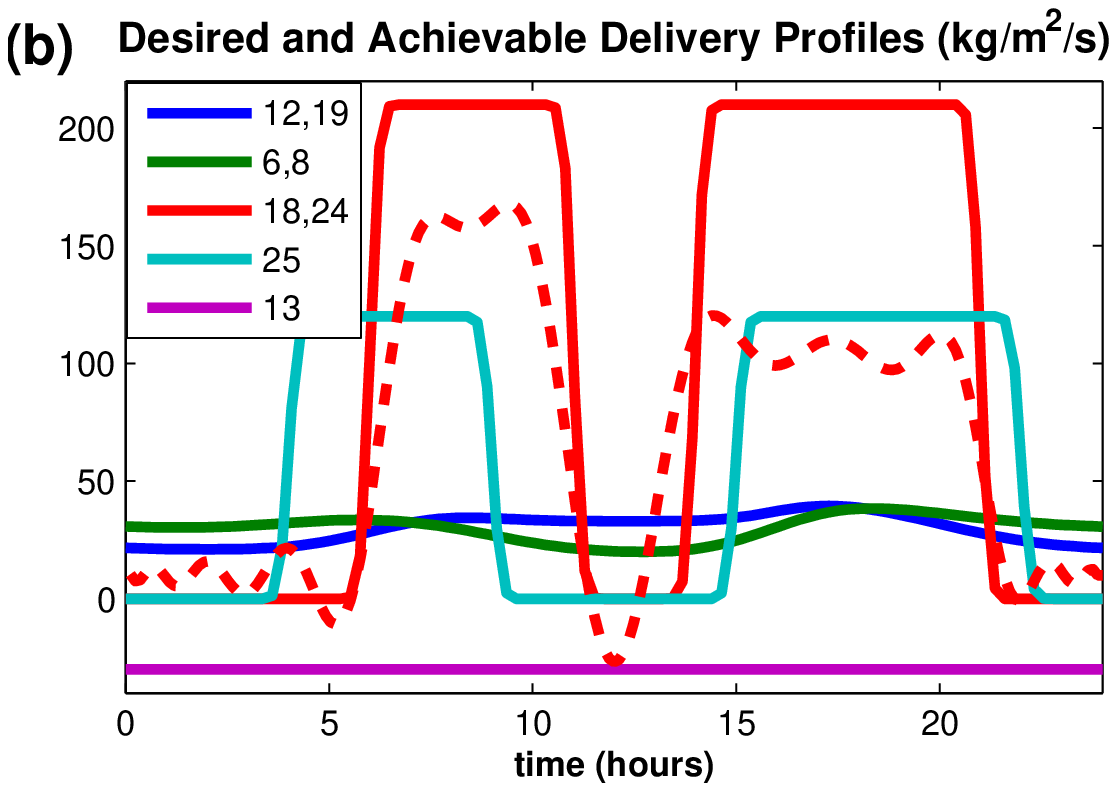} \caption{MLS Optimization. (a) Compression solutions (color), scaled extremal system pressure (thick dashed line), and scaled pressure limits (dotted line); (b) desired (solid) and maximal non-firm (dashed) deliveries.} \label{fig:net4}
\end{figure}

{\it MLS Optimization.} Consider the same scenario as for the ETC problem, except the loads at non-firm customers at nodes $\{18,24\}$ are increased so that there is no feasible solution.  Therefore the MLS objective is used, which adds $(N+1)\times\Sigma=52$ optimization variables.  A priority weighting of $c_{18}(t)=c_{24}(t)=1$ is used.  The desired and maximal non-firm deliveries are shown in Figure \ref{fig:net4}b, and the associated compression solutions are given in Figure \ref{fig:net4}a.  With the right control protocol, the system utilization can be significantly increased even over the delivery profiles in Figure \ref{fig:net3}f and still yield a feasible solution.

\section{Conclusion} \label{sec:conc}

We have developed a new framework for modeling and optimal control of compressible gas flow through pipeline networks with time-varying injections, withdrawals, and control actions of compressors and regulators.  The inclusion of information about transient parameters into a physically representative model and an efficient and tractable optimization scheme are shown to facilitate a dramatic and unprecedented improvement in both the capacity and security of pipeline operations with respect to current industry practice.  The economic transient compression (ETC) formulation provides a cost efficient solution when the desired mass transfer is feasible, and the minimum load shedding (MLS) objective minimizes unfulfilled deliveries if they cannot all be met.  Moreover, the solutions produced by the optimization scheme are validated by direct simulation of the control system model.  In addition, our technique leverages the inherent sparsity of the problem for efficient scaling to larger systems.  Implementation in practice would dramatically increase the effective capacity of gas pipeline systems, and save significant resources now used to provide unutilized margins.

\section*{Acknowledgement} This work was carried out under the auspices of the National Nuclear Security Administration of the U.S. Department of Energy at Los Alamos National Laboratory under Contract No. DE-AC52-06NA25396, and was partially supported by DTRA Basic Research Project \#10027-13399 and by the Advanced Grid Modeling Program in the U.S. Department of Energy Office of Electricity.


\bibliographystyle{unsrt}
\bibliography{gas_master}

\end{document}